\documentclass{gtart_h}  

\def\ifplaintex{\expandafter\ifx\csname documentclass\endcsname\relax}

\ifplaintex 
\else
\fi

\expandafter\ifx\csname beginpicture\endcsname\relax
\expandafter\ifx\csname documentclass\endcsname\relax
\input pictex \else
\input prepictex \input pictex \input postpictex \fi\fi

\def\gt{{\mathsurround=0pt\it $\cal G\mskip-2mu$eometry \&\ 
$\cal T\!\!$opology}}        %  journal title in recommended style

\def\gtp{{\mathsurround=0pt\it $\cal G\mskip-2mu$eometry \&\ 
$\cal T\!\!$opology $\cal P\!$ublications}}  % GT publications

\def\lognumber#1{\def\thelognumber{#1}}
\def\volumenumber#1{\def\thevolumenumber{#1}}
\def\papernumber#1{\def\thepapernumber{#1}}
\def\volumeyear#1{\def\thevolumeyear{#1}}

\def\pagenumbers#1#2{\def\startpage{#1}\def\finishpage{#2}}
\def\published#1{\def\publishdate{#1}}
\def\proposed#1{\def\theproposer{#1}}
\def\seconded#1{\def\theseconders{#1}}
\def\received#1{\def\receiveddate{#1}}
\def\revised#1{\def\reviseddate{#1}}
\def\accepted#1{\def\accepteddate{#1}}

\long\def\asciiabstract#1{\long\def\theasciiabstract{#1}}
\def\asciikeywords#1{\def\theasciikeywords{#1}}

\let\\\par\let\thelognumber\relax
\let\thevolumenumber\relax\let\thepapernumber\relax
\let\thevolumeyear\relax\let\thesamplenumber\relax\let\startpage\relax
\let\finishpage\relax\let\publishdate\relax\let\receiveddate\relax
\let\reviseddate\relax\let\accepteddate\relax\let\theasciititle\relax
\let\theasciiauthors\relax
\let\theasciiabstract\relax\let\theasciikeywords\relax
\let\theasciiemail\relax\let\theshortauthors\relax\let\theshorttitle\relax

\long\def\maketitlep{   % start of definition of \maketitlep

\count0=\startpage

\gt\hfill      %   Journal title (top left) 
\beginpicture
\setcoordinatesystem units <0.33truein, 0.33truein> point at 2.2 0.9
\setplotsymbol ({$\cal G$})
\plotsymbolspacing=9truept
\circulararc 315 degrees from 0 1 center at 0 0
\setplotsymbol ({$\cal T$})
\circulararc 315 degrees from 1 -1 center at 1 0
\endpicture
\break
{\small\ifx\thesamplenumber\relax % sample?  
Volume \else Sample
\fi\thevolumenumber\ (\thevolumeyear)
\startpage--\finishpage\nl
Published: \publishdate}
\vglue 0.5truein plus 0.4fil minus 0.1truein

{\parskip=0pt\leftskip 0pt plus 1fil\def\\{\par\smallskip}{\ifplaintex\large
\else\Large\fi\bf\thetitle}\par\medskip}   

\vglue 0pt plus 0.1fil 

{\parskip=0pt\leftskip 0pt plus 1fil\def\\{\par}{\sc\theauthors}
\par\medskip}

\vglue 0pt plus 0.1fil 

{\small\parskip=0pt\let\newline\\
{\leftskip 0pt plus 1fil\def\\{\par}{\sl\theaddress}\par}
\expandafter\ifx\theemail\relax    % email address?
\relax\else\vglue 5pt plus 0.02fil minus 2pt\def\\{\stdspace{\rm 
and}\stdspace} 
\cl{Email:\stdspace\tt\theemail}\fi
\ifx\theurl\relax                  % URL given?
\relax\else\vglue 5pt plus 0.02fil minus 2pt\def\\{\stdspace{\rm 
and}\stdspace}
\cl{URL:\stdspace\tt\theurl}\fi\par}

\vglue 7pt plus 0.3fil minus 3pt

{\bf Abstract}
\vglue 5pt plus 0.1fil minus 2pt

\theabstract

\vglue 7pt plus 0.3fil minus 3pt

{\bf AMS Classification numbers}\quad Primary:\quad \theprimaryclass

Secondary:\quad \thesecondaryclass

\vglue 5pt plus 0.3fil minus 2pt

{\bf Keywords:}\quad \thekeywords

\vglue 10pt plus 0.5fil minus 5pt

{\small  Proposed: \theproposer\hfill Received: \receiveddate\nl
Seconded: \theseconders\hfill 
\ifx\reviseddate\relax                         % paper revised?
Accepted: \accepteddate                        % no
\else
Revised: \reviseddate                          % yes
\fi}
\eject
}       %  end of definition of \maketitlep

\let\maketitlepage\maketitlep
\let\maketitle\maketitlepage

\font\phead=cmsl9 scaled 950
\font=cmsl9 scaled 1050
\font\pnum=cmbx10 scaled 913
\font\lnum=cmbx10 
\font\pfoot=cmsl9 scaled 950
\font=cmsl9 scaled 1050
\ifplaintex
\headline{\vbox to 0pt{\vskip -4.5mm\line{\small\phead\ifnum
\count0=\startpage ISSN 1364-0380 (on line)
1465-3060 (printed) \hfill {\pnum\folio}\else\ifodd\count0\def\\{ }% 
\ifx\theshorttitle\relax\thetitle\else\theshorttitle\fi\hfill{\pnum\folio}
\else\def\\{ and }{\pnum\folio}\hfill\ifx\theshortauthors\relax\theauthors
\else\theshortauthors\fi\fi\fi}\vss}}
\footline{\vbox to 0pt{\vglue 0mm\line{\small\pfoot\ifnum\count0=\startpage
\copyright\ \gtp\hfill\else
\gt, Volume \thevolumenumber\ (\thevolumeyear)\hfill\fi}\vss
}}
\else
\makeatletter
\def\@oddhead{{\small\lhead\ifnum\count0=\startpage ISSN 1364-0380 (on line)
1465-3060 (printed) \hfill {\lnum\number\count0}\else\ifodd\count0
\def\\{ }\ifx\theshorttitle\relax \thetitle \else\theshorttitle\fi\hfill
{\lnum\number\count0}\else\def\\{ and }{\lnum\number\count0}
\hfill\ifx\theshortauthors\relax 
\theauthors\else\theshortauthors\fi\fi\fi}}\def\@evenhead{\@oddhead}
\def\@oddfoot{\small\lfoot\ifnum\count0=\startpage\copyright\ \gtp\hfill\else
\gt, Volume \thevolumenumber\ (\thevolumeyear)\hfill\fi}
\def\@evenfoot{\@oddfoot}
\makeatother
\fi

\newwrite\gtoutfile
\long\gdef\makeheadfile{  %%% start of definition of \makeheadfile
{\def\\{, }\def\s{ }
\immediate\openout\gtoutfile head.xxx
\immediate\write\gtoutfile{Proxy-for: \ifx\theasciiauthors\relax
\theauthors\else\theasciiauthors\fi\s<\ifx\theasciiemail\relax\theemail\else\theasciiemail\fi>}
\immediate\write\gtoutfile{\noexpand\\}
\immediate\write\gtoutfile{Authors: \ifx\theasciiauthors\relax
\theauthors\else\theasciiauthors\fi}
{\def\\{ }\immediate\write\gtoutfile{Title: \ifx\theasciititle\relax
\thetitle\else\theasciititle\fi}}
\immediate\write\gtoutfile{Subj-class: GT or SG or MG etc}
\immediate\write\gtoutfile{MSC-class: \theprimaryclass\ifx\thesecondaryclass\relax\else, \thesecondaryclass\fi}
\immediate\write\gtoutfile{Journal-ref: Geom. Topol. \thevolumenumber
(\thevolumeyear) \startpage-\finishpage}
\immediate\write\gtoutfile{Comments: Published by Geometry and Topology at}
\immediate\write\gtoutfile{\s\s http://www.maths.warwick.ac.uk/gt/GTVol\thevolumenumber/paper\thepapernumber.abs.html}
\immediate\write\gtoutfile{\noexpand\\}
\immediate\write\gtoutfile{}
\ifx\theasciiabstract\relax
\immediate\write\gtoutfile{\theabstract}\else
\immediate\write\gtoutfile{\theasciiabstract}\fi
\immediate\write\gtoutfile{}
\immediate\write\gtoutfile{\noexpand\\}
\immediate\write\gtoutfile{}
\immediate\closeout\gtoutfile}}  %%% end of definition of \makeheadfile

\def\maketitlepage{\maketitlep\makeheadfile}
\let\maketitle\maketitlepage

\lognumber{366}
\volumenumber{8}\papernumber{17}\volumeyear{2004}
\pagenumbers{675}{699}
\received{30 September 2003}
\accepted{13/Feb/04}
\revised{22 April 2004}
\published{24 April 2004}
\proposed{Robion Kirby}
\seconded{Shigeyuki Morita, Ralph Cohen}

\usepackage{amsmath, amssymb}

\newcommand{\La}{\Lambda} 
\newcommand{\M}{\overline{M}}
\newcommand{\Q}{\mathbb{Q}}
\newcommand{\C}{\mathbb{C}}
\newcommand{\com}{\mathbb{C}}

\newcommand{\Z}{\mathbb{Z}}

\newcommand{\bH}{\mathsf{H}}
\newcommand{\Hc}{\bH^\circ}
\newcommand{\Hb}{\bH^\bullet}
\newcommand{\proj}{\mathbf P}
\newcommand{\cO}{{\mathcal{O}}}
\newcommand{\cE}{\mathcal{E}}
\newcommand{\cP}{\mathcal{P}}

\newcommand{\cR}{\mathcal{R}}
\newcommand{\cS}{\mathcal{S}}
\newcommand{\bA}{\mathsf{A}}
\newcommand{\cs}{{\varsigma}}

\newcommand{\LV}{\Lambda^{\frac\infty2}V}
\newcommand{\gli}{\mathfrak{gl}(\infty)}

\newcommand{\lang}{\left\langle}
\newcommand{\rang}{\right\rangle}
\newcommand{\llang}{\left\langle\left.}
\newcommand{\rrang}{\right\rangle\right.}
\newcommand{\kket}[1]{\right|#1\rang} 
\newcommand{\bbra}[1]{\lang #1 \left|} 
\newcommand{\ket}[1]{\left.\left|#1\right.\rang} 
\newcommand{\bra}[1]{\lang\left. #1 \right|\right.} 

\DeclareMathOperator{\End}{End}

\newcommand{\zz}{{\mathfrak{z}}}
\newcommand{\A}{\mathcal{A}}
\newcommand{\F}{\mathcal{F}}

\newcommand{\vac}{v_\emptyset}

\DeclareMathOperator{\Aut}{Aut}
\DeclareMathOperator{\ev}{ev}
\DeclareMathOperator{\wt}{wt}

\numberwithin{equation}{section}

\newtheorem{Theorem}{Theorem}
\newtheorem{Lemma}{Lemma}

\newtheorem{Proposition}[Lemma]{Proposition}

\theoremstyle{definition}

\def\ssubsection{\medskip\addtocounter{subsection}{1}{\large\bf \thesubsection}\quad}
\def\ssubsubsection{\medskip\addtocounter{subsubsection}{1}{\bf \thesubsubsection}\quad}

\begin{document}
\title{Hodge integrals and invariants of the unknot}
\author{A Okounkov\\R Pandharipande}

\address{Department of Mathematics, Princeton University\\
Princeton, NJ 08544, USA}
\email{okounkov@math.princeton.edu, rahulp@math.princeton.edu}

\begin{abstract}
We prove the Gopakumar--Mari\~no--Vafa formula for special cubic
Hodge integrals. The GMV formula arises from Chern--Simons/string 
duality applied to the unknot in the three sphere. The GMV formula is a 
$q$--analog of the ELSV formula for linear Hodge integrals. We find a system 
of bilinear localization equations relating linear and special cubic Hodge 
integrals. The GMV formula then follows easily from the ELSV formula.
An operator form of the GMV formula is presented in the last
section of the paper.
\end{abstract}
\asciiabstract{%
We prove the Gopakumar-Marino-Vafa formula for special cubic
Hodge integrals. The GMV formula arises from Chern-Simons/string 
duality applied to the unknot in the three sphere. The GMV formula is a 
q-analog of the ELSV formula for linear Hodge integrals. We find a system 
of bilinear localization equations relating linear and special cubic Hodge 
integrals. The GMV formula then follows easily from the ELSV formula.
An operator form of the GMV formula is presented in the last
section of the paper.}

\asciikeywords{Hodge integrals, unknot, Gopakumar-Marino-Vafa formula}


\maketitlepage

\setcounter{section}{-1}
\section{Introduction}
\label{s1}

\ssubsection{}
Let $\M_{g,n}$ be the Deligne--Mumford moduli stack of 
stable curves of genus $g$ with $n$ marked points.
We study here Hodge integrals over $\M_{g,n}$.

Let $L_i$ be the line bundle over $\M_{g,n}$ with 
fiber over the moduli point
$$
[C,p_1,\dots,p_n] \in \M_{g,n}
$$
given by the cotangent space $T^*_{C,p_i}$ of the curve $C$ at $p_i$. 
The $\psi$ classes are the first Chern classes of the
cotangent line bundles,
$$
\psi_i = c_1(L_i) 
$$
 in $H^2(\M_{g,n}, \Q)\,.$

Let $\pi: C \to \overline{M}_{g,n}$ be the universal curve.
Let $\omega_\pi$ be the relative dualizing sheaf.
Let ${\mathbb E}$ be the rank $g$ Hodge bundle on the moduli space
of curves,
$$
{\mathbb E} = \pi_*(\omega_\pi)\,.
$$
The $\lambda$ classes are defined by
$$
\lambda_i = c_i({\mathbb E})
$$
in $H^{2i}(\overline{M}_{g,n}, \Q)$.
The Chern polynomial of the Hodge bundle,
$$
\La(t) = 1 + t \lambda_1 + \dots + t^g \lambda_g,
$$
will appear often in the paper.
 
By definition, the Hodge integrals are the integrals of the
$\psi$ and
$\lambda$ classes over $\overline{M}_{g,n}$.

\ssubsection
The polynomial, 
\begin{equation}
   \label{dHc}
\Hc_g(z_1,\dots,z_n;t_1,\dots,t_s) =  \prod_{i=1}^n z_i \  
\int_{\M_{g,n}} \frac{ \prod_{i=1}^s \Lambda(t_i)}
{ \prod_{i=1}^n (1-z_i \psi_i)}\, ,
\end{equation}
generates  
Hodge integrals on $\M_{g,n}$ with at most
$s$ classes $\lambda_i$. 
In case $s=0$, 
$\Hc_g(z)$ 
generates pure $\psi$ integrals
on $\overline{M}_{g,n}$.
The prefactor $\prod z_i$ 
in the definition is introduced for later
convenience.

The Hodge integral in \eqref{dHc} arises as 
a vertex term in the localization formula for
the Gromov--Witten invariants of an $s$ dimensional target \cite{GP}.

\ssubsection
Mumford's relations may be used to 
reduce Hodge integrals
to $\psi$ integrals containing no
$\lambda$ classes, see \cite{F, FP1,M}. 
However, 
the reduction
method 
often destroys
the rich structure possessed by special 
classes of Hodge integrals.

\ssubsection
Linear Hodge integrals, generated by
$\Hc_g(z;t_1)$, are
connected 
to Hurwitz theory and the related 
combinatorics of symmetric group 
characters. In case
\begin{equation}
\quad z_i \in \mathbb{N}\, , \quad t_1 = -1 \, ,
\label{dc}
\end{equation}
the Ekedahl--Lando--Shapiro--Vainstein (ELSV) formula
expresses $\Hc_g(z;t_1)$ in terms of a Hurwitz 
number \cite{ELSV}. For general
evaluations, there exists an operator
formula expressing    $\Hc_g(z;t_1)$ in terms
of vacuum expectations of 
products of explicit operators in 
Fock space  \cite{P1}. 

The $s=0$ case of pure $\psi$ integrals, studied by Witten
and Kontsevich, is perhaps best 
understood as the $t_1\to 0$ limit 
of the $s=1$ case of linear Hodge integrals, see \cite{OP}. 
In particular, the appearance of random 
matrices should be viewed as the 
continuous limit of random partitions,
which play a central role
in the $s=1$ theory.

\ssubsection
We study here {\em special cubic Hodge integrals}: the Hodge
integrals generated by
$\Hc_g(z;t_1,t_2,t_3)$, 
with parameters $t_i$ 
subject to the 
constraint
\begin{equation}
  \label{CY}
  \boxed{\frac1{t_1}+\frac1{t_2}+\frac1{t_3} = 0} \,.
\end{equation}
In the context of the localization formulas,
the constraint \eqref{CY} is the local 
Calabi--Yau condition. 
The linear Hodge integrals can be recovered
from the special cubic Hodge integrals by a limit:
\begin{equation}
\Hc_g(z;t_1) = \label{lim}
  \lim_{ t_2,t_3\to 0} \Hc_g(z;t_1, t_2,t_3) ,
\end{equation}
where $t_1$ is held fixed and the parameters are
subject to the constraint \eqref{CY}.

We prove two formulas for special cubic Hodge integrals.
The first,  covering the evaluations \eqref{dc},  
is a $q$--analog, or trigonometric 
deformation, of the ELSV formula. The formula is
based upon the conjectural Chern--Simons/string
duality of Gopakumar and Vafa \cite{GoV} applied to 
the case of the unknot in $S^3$ \cite{MV}.
Hodge integrals enter via
a
heuristic localization computation
of S. Katz and C.-C. Liu of the corresponding 
open string Gromov--Witten theory \cite{KL}.
An exposition for mathematicians of the physics behind the 
conjecture 
can be found in \cite{MM}. 
We will call the formula the Gopakumar--Mari\~no--Vafa (GMV)
formula. 

While a precise mathematical definition
of open string Gromov--Witten theory
is lacking at the moment, our results provide
significant theoretical evidence for both 
the general Gopakumar--Vafa conjecture and
the assumptions made in \cite{KL}. 

Our second formula for special cubic Hodge integrals
 is a $q$--analog
of the operator formula of \cite{P1}.  
Our derivation of the operator formula from the GMV formula follows
the steps taken in \cite{P1}.

\ssubsection
In order to state the GMV formula for Hodge 
integrals, we introduce
the generating series 
$$
\Hc(z;t;u) = \sum_{g\ge 0} u^{2g-2} \, \Hc_g(z;t).
$$
We follow the conventions of \cite{P1} regarding the
unstable terms,
$$\Hc_0(z_1;t)= \frac{1}{z_1},$$
$$\Hc_0(z_1,z_2;t)= \frac{z_1z_2}{z_1+z_2}.$$
Let $\Hb(z;t;u)$ denote the disconnected 
$n$--point series.
For example,
$$
\Hb(z_1,z_2;t;u) = \Hc(z_1,z_2;t;u) +
\Hc(z_1;t;u)\, \Hc(z_2;t;u)\,. 
$$
Further discussion of these definitions
can be found in 
\cite{P1}.

The $0$--point function $\Hc(;t;u)$ is {\em not}
included in the disconnected $n$--point series.
For linear Hodge integrals, the $0$--point
function vanishes. 
For cubic Hodge integrals,
the $0$--point function is determined by the results of \cite{FP1}:
\begin{eqnarray*}
\Hc_g(;t_1,t_2,t_3) & =& \sum_{\sigma \in S_3} {t^g_{\sigma(1)} 
t^{g-1}_{\sigma(2)} t^{g-2}_{\sigma(3)}} \int_{\overline{M}_g} \lambda_g \lambda_{g-1}
\lambda_{g-2} +
{t_1^{g-1} t_2^{g-1} t_3^{g-1}} \int_{\overline{M}_g} \lambda_{g-1}^3 \\
& = & \Big( \sum_{\sigma \in S_3}  \frac{t^g_{\sigma(1)} 
t^{g-1}_{\sigma(2)} t^{g-2}_{\sigma(3)}}{2} + {t_1^{g-1} t_2^{g-1} t_3^{g-1}} \Big)
\frac{|B_{2g}|}{2g} \frac{|B_{2g-2}|}{2g-2} \frac{1}{(2g-2)!},
\end{eqnarray*}
for $g\geq 2$.
Here, $B_{m}$ denotes the Bernoulli number.

\ssubsection
A partition $\mu$ is a monotone
sequence
$
 (\mu_1\ge \mu_2 \ge \dots \geq  0)
$
of natural numbers such that $\mu_i=0$ for sufficiently
large $i$.
Let $|\mu|$ denote the size of the partition,
$$|\mu|=\sum \mu_i,$$
and let $\ell(\mu)$ denote the length,
the number of nonzero parts of $\mu$. 
Let  $\Aut \mu$ be the group permuting 
equal parts of $\mu$.

The ELSV formula is 
equivalent to the following equality relating
Hodge integrals to the representation theory of the
symmetric group: 
\begin{equation}
  \label{ELSV}
\left(\prod \frac{\mu_i^{\mu_i}}{\mu_i !} \right)
\Hb(\mu;-1;u) = 
\frac{u^{-|\mu|-\ell(\mu)}}{|\Aut \mu|} 
\sum_{|\lambda|=|\mu|} 
\left(\frac{\dim \lambda}{|\lambda|!} \right)
\, e^{u f_2(\lambda)} \, \chi^\lambda_\mu \,. 
\end{equation}
The summation is over all partitions
$\lambda$ of  size $|\mu|$. 

Here,  
$\dim\lambda$ is the dimension of 
the corresponding representation of the 
symmetric group, $\chi^\lambda_\mu$ is the
irreducible character of the symmetric 
group corresponding to representation $\lambda$
and conjugacy class $\mu$, and $f_2(\lambda)$ 
is the central character 
$$
f_2(\lambda) = \binom{|\lambda|}{2} 
\, \frac{ \chi^\lambda_{(2,1,\dots,1)}}{  \dim \lambda}
$$ 
of a transposition 
in the representation $\lambda$. Explicitly, the 
central character is given by the following formula 
$$
f_2(\lambda) = \frac12 \sum_{i\geq 1} 
\left[(\lambda_i-i+\frac{1}{2})^2 - (-i+\frac{1}{2})^2\right] \,.
$$

\ssubsection
The GMV formula is a $q$--deformation of the
ELSV formula \eqref{ELSV}. In  fact, the only term
of the right side of formula \eqref{ELSV} which is 
deformed is the dimension $\dim \lambda$. The dimension
$\dim \lambda$ is 
replaced by the $q$--dimension of 
$\lambda$,  a well-known notion in the 
theory of quantum groups and the corresponding theory of 
knot invariants, see for example \cite{KA}. 

The $q$--dimension of $\lambda$ is the rational function of the
variable $q^\frac12$ defined
by
\begin{equation}
\dim_q \lambda = q^{-\frac12 f_2(\lambda) - \frac12 |\lambda|} \,
s_\lambda(1,q^{-1},q^{-2},q^{-3},\dots)\,,
\label{dqdim}
\end{equation}
where $s_\lambda$ is the Schur function.
Alternatively, the $q$--dimension may be defined by
%This is a rational function 
%of variable $q^{1/2}$ which is explicitly given by
%the following formula
% 
\begin{equation}
    \label{qhook}
\frac{ \dim_q \lambda}{|\lambda|!} 
= \prod_{\square \in \lambda} \frac{1}
{q^{h(\square)/2}-q^{-h(\square)/2}} \,, 
\end{equation}
where the product is over all squares $\square$ in 
the diagram of $\lambda$ and $h(\square)$ is the
corresponding hook-length. The standard hook-length formula
for $\dim \lambda$ 
arises as the coefficient of the leading term in \eqref{qhook}
as
$q\to 1$.

\ssubsection
The first result of the paper is the Gopakumar--Mari\~no--Vafa formula.
The proof is presented in Section \ref{prr}.

\begin{Theorem}[Gopakumar--Mari\~no--Vafa formula]
For any number $a$ and any partition $\mu$ we have 
  \begin{multline}\label{GV}
\left(\prod 
\binom{(a+1)\mu_i}{\mu_i} \right)
\, \Hb\left(\mu;-1,-\frac1a,\frac1{a+1};\sqrt{a(a+1)}\, u
\right) = \\
\frac{(au)^{-\ell(\mu)}}  {|\Aut \mu|} \sum_{|\lambda|=|\mu|}
\left(\frac{ \dim_{q} \lambda}{|\lambda|!} \right)\,
e^{(a+\frac12) u f_2(\lambda)} \, \chi^\lambda_\mu \,,
  \end{multline}
where $q=e^u$.   
\end{Theorem}

The $q$--dimension of $\lambda$ may be rewritten after
the substitution $q=e^u$ as 
\begin{equation}
  \label{qhu}
  \frac{ \dim_{\,e^u}\! \lambda}{|\lambda|!} = 
\prod_{\square \in \lambda} \,\, \left(
2 \sinh \frac{u h(\square)}2 \right)^{-1} \,.
\end{equation}
A straightforward analysis shows the
GMV formula specializes to the ELSV formula
 as $a\to\infty$ and
$u\to 0$ while keeping the product $au$
constant.

\ssubsection
Following the treatment of linear Hodge integrals 
in \cite{P1}, the GMV formula can be expressed in
operator form, 
\begin{equation}\label{Hoper}
    \Hb\left(z;-1,-\frac1a,\frac1{a+1}; \sqrt{a(a+1)}\, u
\right) = 
\lang
\prod \bA(z_i;a) \rang  \,.
\end{equation}
Here, 
$\bA(z_i;a)$ is an explicit operator in a 
fermionic Fock space, and the
angle brackets denote the vacuum expectation.
The operator formula, 
Theorem \ref{thm2}, is fully stated and proven
in Section \ref{s2}.

%The validity of  
%formula 
%\eqref{Hoper} 
%for all values 
%of the arguments $z$ is a significant advantage.
The
operators $\bA(z_i;a)$ are $q$--analogs of the operators
studied in \cite{P1} in connection with linear 
Hodge integrals. The operator formula for linear Hodge 
integrals played a fundamental role in the study of
the equivariant Gromov--Witten  
of $\proj^1$ in \cite{P1}. 
Similarly, the operator formula for special cubic Hodge integrals
is applicable to the study of several local Calabi--Yau
geometries.

\subsection{Acknowledgements}
The results of the paper were presented by A.O. in
a workshop on {\em Moduli, Lie theory, and Interactions with
Physics} sponsored by the Istituto Nazionale di Alta Matematica in
Rome in the spring of 2003. Special thanks are due to the
organizers for their hospitality.
A 
different proof the GMV formula can be found in
\cite{LLZ1}, \cite{LLZ2}, see also \cite{LLZ3}.  

A.O. was partially supported by 
DMS-0096246 and fellowships from the Sloan and Packard foundations.
R.P. was partially supported by 
DMS-0236984 and fellowships from the Sloan and Packard foundations.

\section{Proof of the GMV formula}
\label{prr}

\subsection{Strategy}

We find a system of bilinear localization equations which
relates linear Hodge integrals to special cubic
Hodge integrals with the following two properties:
\begin{enumerate}
\item[(i)] given the linear Hodge integrals, the system has a
           unique 
solution for the complete set of special cubic Hodge integrals,
\item[(ii)] the ELSV and GMV formulas satisfy the system.
\end{enumerate}
Localization equations of a similar form 
were used successfully in \cite{FP2}. Our strategy here
was motivated by \cite{FP2}.

In fact, bilinear localization equations uniquely constrain
{\em all} cubic Hodge integrals. A further study of
cubic Hodge integrals will be presented in a future paper.

\subsection{Integrals}

\ssubsubsection
Our localization relations are obtained from the analysis of
integrals over the moduli space of stable maps to $\proj^1$.

Let  
$\M_{g,n}(\proj^1,d)$ be the moduli space of degree $d>0$ stable maps
from genus $g$, $n$--pointed {\em connected} curves to $\proj^1$.
The virtual dimension of the moduli space $\M_{g,n}(\proj^1,d)$ is
$2g+2d-2$.
Let
$$
\ev_i: \M_{g,n}(\proj^1,d) \to \proj^1
$$
denote the evaluation map at the $i$th marked point, and 
let
$$\pi: C \rightarrow \M_{g,n}(\proj^1,d),$$
$$f: C \rightarrow \proj^1$$
denote the universal curve and universal map.
The bundles
$${\mathbb A} = R^1\pi_*f^* \cO, \ \ 
{\mathbb B} = R^1 \pi_*f^* \cO(-1),$$
will play an important role. 
The ranks of ${\mathbb A}$ and ${\mathbb B}$
are $g$ and $g+d-1$ respectively.

Let 
$
\nu=(\nu_1,\dots,\nu_n) 
$
be a partition such that
$
|\nu| \le d-1 \,.
$
Here, we allow $\nu_i=0$.
Consider the  integral
\begin{equation}
    \label{defI}
I_{g,d}(\nu) = \int_{[\M_{g,n}(\proj^1,d)]^{vir}}
\ev_1^*(\omega)^{d-1-|\nu|} \, \, 
\prod_{i=1}^{n} \psi_i^{\nu_i} \ev_i^*(\omega)\,
\, c_{top}({\mathbb A}) \, 
c_{top}({\mathbb B})\, ,
\end{equation}
where $\omega\in H^2(\proj^1, \com)$ is the Poincar\'e dual of
the point class. The dimension of the integrand equals the
the virtual dimension of  $\M_{g,n}(\proj^1,d)]$.

The integrand in \eqref{defI} involves 
the class $\ev^*_1(\omega)^{d-|\nu|}$. Since 
$$
\omega^2=0,
$$
$I_{g,d}(\nu)$  vanishes if $|\nu|<d-1$. 

If $|\nu|=d-1$, the above vanishing does not apply.
Let
$$
I_d(\nu;u) = \sum_{g\ge 0} u^{2g-2} \, I_{g,d}(\nu) \,.
$$
The nonvanishing 
values of the series $I_{d}(\nu;u)$ are determined
by the following result proven in Section \ref{sp1}.

\begin{Proposition}\label{p1}
If $|\nu|=d-1$, then
$$
I_d(\nu;u) = (-1)^{d+1} \frac{d^{n-2}}{\prod \nu_i!} 
\, \frac{1}{2 u\sin \frac{du}{2}} \,.
$$
\end{Proposition}

\subsection{Virtual localization}

\ssubsubsection
Let $\com^*$ act diagonally on a two dimensional vector
space $V$ via the trivial and
standard representations,
\begin{equation}
\label{repp}
\xi\cdot (v_1,v_2) = ( v_1, 
\xi \cdot v_2).
\end{equation}
Let $\proj^1=\proj(V)$.
Let $0,\infty $ be the fixed points $[1,0], [0,1]$ of the corresponding
$\com^*$--action on $\proj(V)$.

An equivariant lifting  of $\com^*$ to a line bundle $L$
over 
$\proj(V)$ is uniquely determined by
the fiber
representations $L_{0}$ and $L_{\infty}$ at the fixed points.
%We denote 1--dimensional representations of $\com^*$ by 
%their integer  weight.
We represent the data of the $\com^*$ lifting
to $L$ by the integral weights $[l_0,l_\infty]$ 
of the representations $L_0$ and $L_\infty$.  
The canonical lifting of $\com^*$ to the
tangent bundle $T_{\proj^1}$ has weights $[1,-1]$.

Let $a$ be an integer. Equivariant lifts of $\com^*$ to
the line bundles $\cO$ and $\cO(-1)$ are given
by the weights
\begin{equation}
\label{wfw}
[-a,-a], \ \ [a,a+1],
\end{equation}
respectively.

The Poincar\'e dual of a point, $\omega \in H^2(\proj^1,\com)$, can lifted
to the equivariant cohomology of $\proj^1$ by selecting the class of either
fixed point $0,\infty \in \proj^1$. We will first define the lift
by
\begin{equation} \label{wwffww}
\omega = [0] \in H^2_{\com^*}(\proj^1,\com).
\end{equation}
The second choice will be used in Section \ref{sp1}.

The representation (\ref{repp}) canonically
induces a $\com^*$--action on the moduli space 
$\overline{M}_{g,n}(\proj^1,d)$  by translation of
the map. The $\com^*$--action lifts canonically
to the cotangent lines of the moduli space of
maps. 
A localization formula for the  virtual fundamental class is proven
in \cite{GP}.

A $\com^*$--equivariant lift of the integrand of $I_{g,d}(\nu)$ is
defined for each parameter $a$
by the lifts (\ref{wfw}) and (\ref{wwffww}) .
The virtual localization formula provides an evaluation of
$I_{g,d}(\nu)$ in terms of Hodge integrals.

\ssubsubsection
The localization analysis of $I_{g,d}(\nu)$ is uniform for
$a \neq -1,0$.
The localized equivariant
integral,
$\tilde{I}_{g,d}(\nu)$, defined by
$$
\int_{[\M_{g,n}(\proj^1,d)]^{vir}}\!
\ev_1^*(\omega)^{d-1-|\nu|} \, 
\prod_{i=1}^{n} \psi_i^{\nu_i} \ev_i^*(\omega)\,
\, \frac{c_{top}(R^1\pi_*f^* \cO)}{c_{top} (R^0\pi_*f^* \cO)} \, 
\frac{c_{top}(R^1\pi_*f^*\cO(-1))}{
c_{top}(R^0\pi_*f^*\cO(-1))}
\, ,$$
is more convenient to study in case $a$ is nonzero.
Since,
$$ -\frac{1}{a} {I}_{g,d}(\nu)= \tilde{I}_{g,d}(\nu),$$
the difference is quite minor.

\ssubsection
The virtual localization formula expresses 
equivariant integrals over the moduli space 
$\overline{M}_{g,d}(\proj^1,d)$ as a sum over
localization graphs $\Gamma$. 
The graphs $\Gamma$ corresponds bijectively to
components of the locus of the 
$\C^*$--fixed points in $\M_{g,n}(\proj^1,d)$.
A complete discussion of
the graph structure of the virtual localization formula
for $\overline{M}_{g,n}(\proj^1,d)$ can be found in 
\cite{OP}.

Let $[f]\in \overline{M}_{g,n}(\proj^1,d)$ represent
a generic point of a component of the $\com^*$--fixed locus,
$$
[f:(C,p_1,\dots,p_n) \to \proj^1].
$$
The components of $C$ are either collapsed by $f$ to 
$\{0,\infty\}\in\proj^1$ or are unmarked
rational Galois covers of $\proj^1$ fully ramified over
$0$ and $\infty$.

The vertex set $V$ of $\Gamma$ corresponds to the connected components
of the set,
$$f^{-1}( \{ 0, \infty\}).$$
Excepting degenerate issues, the vertices 
correspond to
the collapsed components of $C$.
The vertices
carry
a genus $g_v$, a marking set, 
and an assignment to $0$ or $\infty$ in $\proj^1$.

By our choice of 
the equivariant lift of the class $\omega$, only vertices
lying over $0\in\proj^1$ are allowed to carry markings.
Graphs $\Gamma$ with markings on vertices lying over $\infty$
do not contribute to the localization calculation of $\tilde{I}_{g,d}(\nu)$. 
 
The edge set $E$ of $\Gamma$ corresponds
to the non-collapsed components $C_i$ of $C$. The edges 
carry the degree $\mu_i$ of the 
map $f\big|_{C_i}$. The edge degrees  $\mu_i$ form
a partition of the number $d$. 

The localization graphs $\Gamma$ for $\overline{M}_{g,n}(\proj^1,d)$ must
be connected and satisfy the global genus 
condition
$$
\sum_{\textup{vertices $v$}} g_v + h^1(\Gamma)= g \,,
$$
where $h^1(\Gamma)$ is the first Betti number of $\Gamma$ 
$$
h^1(\Gamma) = 1-|V| + |E|. 
$$

\ssubsubsection
Let $\Gamma$ be a localization graph for the moduli space
$\overline{M}_{g,n}(\proj^1,d)$. 
The localization contribution of the graph $\Gamma$ to $\tilde{I}_{g,d}(\nu)$
factors into vertex and edge 
contributions.

\begin{enumerate}
\item[$\bullet$]
Let $v$ be a vertex lying over $0\in\proj^1$
carrying genus $g_v$ and 
$s$ marked points $p_1,\dots,p_s$.
Let $\mu_1,\dots, \mu_r$ be the degrees associated
to the edges incident to $v$. 
The contribution of $v$
 to the localization formula for $\tilde{I}_{g,d}(\nu)$ 
is: 
\begin{eqnarray*}
\text{Cont}(v) & = & (-1)^{g_v-1} a^{2g_v-2} \int_{\M_{g_v,s+r}} 
\frac{\Lambda(-1) \, \Lambda(1/a) \, \Lambda(-1/a)\,
\prod_{1}^s \psi_i^{\nu_i}}
{\prod_{i=1}^r \left(1\big/{\mu_i} - \psi_{s+i}\right)} \\ &= &
\left[u^{2g_v-2} \prod z_i^{\nu_i+1} \right] \,
\Hc(z_1,\dots,z_s,\mu_1,\dots,\mu_r;-1;iau) \,.
\end{eqnarray*}
\end{enumerate}
Here, $\left[u^{2g_v-2} \prod z_i^{\nu_i+1} \right] \,
\Hc$ denotes the coefficient of the monomial in the function $\Hc$.
The second equality above follows from 
 the relation,
\begin{equation}
  \label{Mumf}
  \Lambda(t) \, \Lambda(-t) = 1 \,,
\end{equation}
proven by Mumford \cite{M}.
\begin{enumerate}
\item[$\bullet$] Let $v$ be a vertex lying over $\infty$
carrying genus $g_v$. Let  
$\mu_1,\dots, \mu_r$ be the degrees associated
to the edges incident to $v$. The contribution of $v$
to the 
localization formula for $\tilde{I}_{g,d}(\nu)$ is:
\begin{equation*}
 \text{Cont}(v)=   
 \left[u^{2g_v-2}\right] \, \Hc(\mu_1,\dots, \mu_r; -1,-\frac1a,\frac1{a+1}; 
i \sqrt{a(a+1)} u)
\end{equation*}
\end{enumerate}
As previously noted, no markings are allowed on vertices over $\infty$.
\begin{enumerate}
\item[$\bullet$]
Let $e$ be an edge of degree $\mu$. The contribution of $e$ to the 
localization formula for $\tilde{I}_{g,d}(\nu)$ is:
\begin{equation*}
 \text{Cont}(e)= (-1)^{\mu} \, a^2 \, \frac{\mu^{\mu}}{\mu!} \, 
\binom{(a+1)\mu}{a \mu} \,. 
\end{equation*}
\end{enumerate}
The vertex and edge contributions are obtained directly from the
virtual localization formula, see \cite{GP, OP}.

Let $G_{g,n}(\proj^1,d)$ denote the set of localization graphs for the
moduli space $\overline{M}_{g,n}(\proj^1,d)$. 
The localization formula for $\tilde{I}_{g,d}(\nu)$ is:
$$\tilde{I}_{g,d}(\nu) = \sum_{\Gamma} \frac{1}{|\text{Aut}\ \Gamma|}
\prod_{v\in V} \text{Cont}(v) \prod_{e\in E} \text{Cont}(e).$$
We have proven the following result.

\begin{Lemma} \label{qqp} For $a\neq -1,0$,
$$-\frac{1}{a} I_{g,d}(\nu) =  \sum_{\Gamma\in G_{g,n}(\proj^1,d)}
 \frac{1}{|\text{\em Aut}\ \Gamma|}
\prod_{v\in V} \text{\em Cont}(v) \prod_{e\in E} \text{\em Cont}(e).$$
\end{Lemma}

\label{jjj}

\ssubsubsection
Define the disconnected bilinear Hodge integral function 
$Z^\bullet_d(\nu;u)$ by
the following formula: 
\begin{multline}
  \label{Zloc}
  Z^\bullet_d(\nu;u) = 
(-1)^d \left[\prod z_i^{\nu_i+1}\right] 
\sum_{|\mu|=d} \frac{(au)^{2\ell(\mu)}}{\zz(\mu)} 
\, \prod \frac{\mu_i^{\mu_i}}{\mu_i!} \, 
\binom{(a+1)\mu_i}{a \mu_i} \times \\
\Hb(z,\mu;-1;iau) \, 
\Hb\left(\mu; -1,-\frac1a,\frac1{a+1}; 
i \sqrt{a(a+1)} u\right)\, ,
\end{multline}
where
$$
\zz(\mu)= \left|\Aut(\mu)\right| \prod_{i=1}^{\ell(\mu)} \mu_i \,.
$$
Let $Z^\circ_d(\nu;u)$ denote the connected part of $Z^\bullet_d(\nu;u)$.

Lemma \ref{qqp} together with the formula for the vertex and edge contributions
in Section \ref{jjj} yields the following result.

\begin{Lemma} \label{pqq} For $a\neq -1,0$,
\begin{equation}
-\frac{1}{a} I_d(\nu;u) = Z^\circ_d(\nu;u) \ . 
\label{FZ}
\end{equation}
\end{Lemma}

\subsection{Proof of Proposition \ref{p1}}\label{sp1}

\ssubsubsection
We will evaluate the integral
$I_{g,d}(\nu)$ for $|\nu|=d-1$
by virtual localization. A {\em new} lift of
the $\com^*$-action to the integrand will be used to
evaluate the localization graph sum.

A lift of the $\com^*$--action to the
integrand is chosen as follows.
The lift of the $\com^*$--action to the bundles ${\mathbb A}$ 
and ${\mathbb B}$ is  
defined by the parameter value $a=0$ in (\ref{wfw}).
The class $\omega$ is lifted by 
$$\omega=[\infty] \in H^2_{\com^*}(\proj^1,\com).$$

Let $\Gamma$ be a localization graph with nonvanishing contribution
to the integral $I_{g,d}(\nu)$ with the specified lift:
\begin{enumerate}
\item[(i)] The weight 0 linearization of $\cO(-1)$ over 0 implies
            each vertex of $\Gamma$ of $0$ is of valence $1$.
\item[(ii)] Each 
vertex $v$ over $0$ carries the class $c_{g(v)}({\mathbb E})^2$
            obtained weight zero linearizations of $\cO$ and $\cO(-1)$ 
            over 0. Since 
$$c_{g(v)}({\mathbb E})^2 =0$$
for $g(v)>0$, all vertices over $0$ must have genus 0.
\item[(iii)] Since $\Gamma$ is connected, there is a unique vertex
            $v_\infty$ over $\infty$. The vertex $v_\infty$ carries
             the full genus $g$.
\item[(iv)] All the markings of $\Gamma$ lie on $v_\infty$.
\end{enumerate}
The contributing graphs $\Gamma$ are therefore indexed uniquely by the degree
partition $\mu$ specified by the edges. 

The
localization graph sum directly yields a formula for $I_{g,d}(\nu)$
in terms of a sum over partitions.

\begin{Lemma}
\label{ggss} If $|\nu|=d-1$, then
\begin{equation}
  \label{Igsum}
  I_{g,d}(\nu) = \sum_{|\mu|=d} \frac{(-1)^{\ell(\mu)+1}}{|\Aut\mu|} \, 
\prod \frac{\mu_i^{\mu_i-1}}{\mu_i!} \, \int_{\M_{g,n+\ell(\mu)}}
\lambda_g \, \frac{\prod_{i=1}^n \psi_i^{\nu_i}}
{\prod_{i=1}^{\ell(\mu)} (1-\mu_i \psi_{i+n})}\,.
\end{equation}
\end{Lemma}

\ssubsubsection
The value of the $\lambda_g$--integral on  right side of 
\eqref{Igsum} can be computed using the following formula \cite{FP2},
\begin{equation}
  \label{lamg}
  \int_{\M_{g,m}} \lambda_g \prod_{i=1}^m \psi_i^{\gamma_i} =
\binom{2g-3+m}{\gamma_1,\dots,\gamma_m} \, [u^{2g-2}] \frac{1}{2u \sin u/2} \,.
\end{equation}
Since $|\nu|=d-1$ and $|\mu|=d$, we obtain
\begin{multline}
  \label{lamg2}
   \int_{\M_{g,n+\ell(\mu)}}
\lambda_g \, \frac{\prod_{i=1}^n \psi_i^{\nu_i}}
{\prod_{i=1}^{\ell(\mu)} (1-\mu_i \psi_{i+n})} = \\
d^{2g-2+n-d+\ell(\mu)} 
\binom{d-1}{\nu_1,\dots,\nu_n} \binom{2g-3+n+\ell(\mu)}{d-1} 
\, [u^{2g-2}] \frac{1}{2u \sin u/2} \,. 
\end{multline}

\ssubsubsection
\begin{Lemma} \label{l1}
Let $t$ be a variable, and let
 $k\geq 0$ be an integer. Then,
  \begin{equation}
    \label{ident1}
    \sum_{|\mu|=d} \frac{(-t)^{\ell(\mu)}}{|\Aut \mu|} \,
 \binom{\ell(\mu)}{k} \,   
\prod \frac{\mu_i^{\mu_i-1}}{\mu_i !}  = 
(-1)^{k} \, \frac{t^k (k-t) (-t+d)^{d-k-1}}{k! (d-k)!} \,.
  \end{equation}
Evaluation at $t=d$ yields
\begin{equation}
\begin{cases}
  \dfrac{(-d)^d}{d!}\,,& k= d-1 \,,\\
0\,, & {k=0, 1, \ldots,d-2} \,. 
\end{cases}
\label{resl1}
\end{equation}
\end{Lemma}

%Note that because of the vanishing in \eqref{resl1}, the precise form of the
%lower terms in \eqref{lowt} is immaterial. 

\begin{proof}
  Consider the following function 
  \begin{equation}
    \label{defT}
    T(x) = \sum_{n>0} \frac{n^{n-1}}{n!} \, x^n \,,
  \end{equation}
which enumerates rooted trees with $n$ vertices 
and solves the functional equation
\begin{equation}
x \exp(T(x)) = T(x) \,,\label{feT}
\end{equation}
as shown, in particular, in \cite{Eul}. More
generally, 
%see \cite{PS}, III.~Abschn., Kap.~5, Aufgabe 210,
%that 
%
\begin{equation}
  \label{feT2}
  \exp(t\, T(x)) = \sum_{n\geq 0} \frac{t(t+n)^{n-1}}{n!} \, x^n \, ,
\end{equation}
see \cite{PS}.
The left  side of \eqref{ident1} equals the coefficient
of $x^d$ in the expansion of
\begin{align*}
  \sum_{l\geq k} \frac{(-t\, T(x))^l}{k!(l-k)!} &= \frac{(-1)^k}{k!} t^k T(x)^k
  \exp(-t\, T(x)) \\
& =  \frac{(-1)^k}{k!} t^k \sum_{n} \frac{(k-t)(k-t+n)^{n-1}}{n!} \, x^{n+k}\,,
\end{align*}
where the second equality follows from  \eqref{feT} and \eqref{feT2}. 
\end{proof}

\ssubsubsection
We view the binomial coefficient,
\begin{equation}
\binom{2g-3+n+\ell(\mu)}{d-1} \label{lowt}\ ,
\end{equation}
as a polynomial in the variable $\ell(\mu)$.
The polynomial \eqref{lowt} agrees in highest order with the
polynomial
$$\binom{\ell(\mu)}{d-1}$$
The lower order terms can be expressed as
a linear combination of the polynomials 
$$
\binom{\ell(\mu)}{k}\,,\quad k=0,1,\dots,d-2 \,.
$$

After substituting  the evaluation \eqref{lamg2} into 
the partition sum \eqref{Igsum} and applying Lemma \ref{l1},
 we obtain
\begin{equation}
  \label{Ig2}
   I_{g,d}(\nu) = (-1)^{d+1} \frac{d^{2g-2+n-1}}{\prod \nu_i!} \, 
   [u^{2g-2}] \frac{1}{2u \sin u/2} \,. 
\end{equation}
Hence,
$$
I_d(\nu;u) = (-1)^{d+1} \,\frac{d^{n-2}}{\prod \nu_i!} \, 
   \frac{1}{2u \, \sin \frac{du}2} \,,
$$
concluding the proof of Proposition \ref{p1}. 
\qed

\subsection{Bilinear relations}
Let $a\neq -1,0$ be a fixed integer.
We have found homogeneous and
inhomogeneous bilinear localization equations
relating linear Hodge integrals to the set of special cubic Hodge integrals,
$$
\Hc(\mu_1,\dots, \mu_r; -1,-\frac1a,\frac1{a+1}; 
i \sqrt{a(a+1)} u),$$
indexed by partitions $\mu$.

Let $|\nu|<d-1$. Since $I_d(\nu;u)$ vanishes,
the homogeneous bilinear equation,
\begin{equation}
\label{tttr}
Z^\circ_{d}(\nu;u)=0,
\end{equation}
is obtained from Lemma \ref{pqq}.

Let $|\nu|=d-1$. By Lemma \ref{pqq} and Proposition \ref{p1},
we obtain the inhomogeneous bilinear equation,
 \begin{equation}
\label{tttrr}
Z^\circ_{d}(\nu;u)=-\frac{1}{a}  (-1)^{d+1} \,\frac{d^{n-2}}{\prod \nu_i!} \, 
   \frac{1}{2u \, \sin \frac{du}2},
\end{equation}

\begin{Lemma}
The  bilinear equations \eqref{tttr} and \eqref{tttrr} 
uniquely determine the special cubic Hodge integrals 
\begin{equation}
\label{kkkk}
[u^{2g-2}] \ \Hc(\mu_1,\dots, \mu_r; -1,-\frac1a,\frac1{a+1}; 
i \sqrt{a(a+1)} u)
\end{equation}
 from
 linear Hodge integrals.
\end{Lemma}

\begin{proof}
We proceed by induction on the genus $g$ and the degree $|\mu|$.
The base case and the induction step are proven
simultaneously.

Assume the Lemma is true for all $g'<g$ and all partitions
$\mu'$ of $d'<d$.
Consider first the genus $g$ homogeneous equations,
\begin{equation}
\label{pssp}
[u^{2g-2}] \ Z^\circ_{d}(\nu;u)=0,
\end{equation}
for $|\nu|<d-1$.
We need only consider
the {\em principal terms} corresponding localization graphs
with a single genus $g$ vertex over $\infty$ incident to
all edges. All non-principal terms are
determined by the induction hypothesis.

As the partition $\nu$ varies, we obtain linear equations for 
the scaled vertex integrals of the principal terms,
\begin{equation*}
[u^{2g-2}] \ \frac{(-1)^{d+\ell(\mu)}}  
 {\zz(\mu)} 
\, \prod \frac{\mu_i^{\mu_i}}{\mu_i!} \, 
\binom{(a+1)\mu_i}{a \mu_i} 
%\Hb(z,\mu;-1;iau) \, 
%(iau)^{-2\ell(\mu)}
\Hc\left(\mu; -1,-\frac1a,\frac1{a+1}; 
i \sqrt{a(a+1)} u\right)\, .
\end{equation*}
We view the above scaled vertex integrals as
a set of variables indexed by partitions $\mu$ of $d$.
Since lower terms appear, the equation are {\em not}
homogeneous in the scaled vertex integrals of the principal terms.

We now consider the homogeneous localization equations obtained in case
the parameter $a$ is set to $0$ in the integrand of $I_{g,d}(\nu)$ with
 lift 
$$\omega = [0] \in H^2_{\com^*}(\proj^1,\com).$$
The $a=0$ equations were studied in \cite{FP2} to calculate
$\lambda_g$ integrals.  
The coefficients of the scaled vertex integrals in the linear equations
obtained from \eqref{pssp} for $a\neq -1,0$
are identical to
the coefficients of the scaled $\lambda_g$ integrals,
$$\ \frac{(-1)^{d+\ell(\mu)}}  
 {\zz(\mu)} 
\, \prod \frac{\mu_i^{\mu_i}}{\mu_i!} 
\prod \mu_i
\int_{\overline{M}_{g,\ell(\mu)}}
 \frac{\lambda_g}{\prod (1-\mu_i \psi_i)},$$
in
the linear equations considered in 
\cite{FP2}.
By the main result of \cite{FP2}, the system of linear equations
has a rank 1 solution space for the
scaled vertex integrals.

Next, we study the genus $g$ inhomogeneous equations for $a\neq -1,0$,
\begin{equation}
\label{yppy}
[u^{2g-2}] \ Z^\circ_{d}(\nu;u)= 
[u^{2g-2}]\ \Big( -\frac{1}{a}  (-1)^{d+1} \,\frac{d^{n-2}}{\prod \nu_i!} \, 
   \frac{1}{2u \, \sin \frac{du}2} \Big),
\end{equation}
where $|\nu|=d-1$.
Again we consider the linear equations for the
scaled vertex integrals of the principal terms.

The coefficients of the linear equation for the scaled vertex
integrals obtained from \eqref{yppy} for $a\neq -1,0$ match the corresponding
coefficients of scaled $\lambda_g$ integrals in 
the $a=0$ calculation of 
$I_{g,d(\nu)}$ with lift 
$$\omega = [0] \in H^2_{\com^*}(\proj^1,\com).$$
The $a=0$ calculation {\em consists only of principal terms}.
The scaled $\lambda_g$ integrals in the $a=0$ are {\em not}
annihilated by the $|\nu|=d-1$ equations since the
coefficient
$$[u^{2g-2}]\ \Big( -\frac{1}{a}  (-1)^{d+1} \,\frac{d^{n-2}}{\prod \nu_i!} \, 
   \frac{1}{2u \, \sin \frac{du}2} \Big)$$
is proportional a nonvanishing Bernoulli number, see \cite{FP1}.

The linear equations for the scaled vertex integrals
of the principal terms obtained from \eqref{yppy} for $a\neq -1,0$ 
are therefore
{\em not} dependent upon the linear equations obtained
from \eqref{pssp} for $a\neq -1,0$. 
Therefore, the full set of linear equations determines
the genus $g$ degree $d$ special cubic Hodge integrals
\eqref{kkkk}. 

\end{proof}

For fixed genus, the $a$ dependence of the special cubic
Hodge integral
\begin{equation}
\label{kkkkk}
[u^{2g-2}] \ \Hc(\mu_1,\dots, \mu_r; -1,-\frac1a,\frac1{a+1}; 
i \sqrt{a(a+1)} u)
\end{equation}
is a rational function.  Rational functions are specified
by values taken on integers not equal to $-1,0$.

To prove the GMV formula for special cubic Hodge integrals,
we need only show the ELSV and GMV formulas satisfy
the homogeneous and inhomogeneous localization
equations for all parameters $a\neq -1, 0$.

\subsection{Operator formalism}

\ssubsubsection
We will prove the ELSV and GMV formulas satisfy the
bilinear localization equations by a calculation in
the infinite wedge representation $\LV$.

We refer the reader to \cite{P1} for a discussion of
the vector space $\LV$. We will use several 
fundamental operators on $\LV$
defined in \cite{P1}.

\ssubsubsection
The first step is to rewrite $Z^\bullet_d(\nu,u)$
defined by formula \eqref{Zloc} 
as a vacuum expectation in $\LV$. 

An operator formula for Hodge integrals of 
the form $\Hb(z;-1;u)$ was obtained in 
\cite{P1}. Applied to the present setting, we find:
\begin{multline}
  \label{oper} 
\Hb(z_1,\dots,z_n,\mu;-1;iau) =\\
(iau)^{-n-d-\ell(\mu)} \left(\prod \frac{\mu_i!}{\mu_i^{\mu_i}}
\right) \llang \prod \A(z_i,iau z_i) \, e^{\alpha_1}  \, e^{iau \F_2} 
\, \kket{\mu}\, ,
\end{multline}
where
\begin{equation}
\ket{\mu}=\prod \alpha_{-\mu_i} \, \vac \,.
\label{mu}
\end{equation} 

The GMV formula \eqref{GV} can also be recast in an 
operator form. The Schur function entering the 
definition of the $q$--dimension can be 
written as the following matrix element: 
\begin{equation}
  \label{SchuG}
  (\Gamma_+(u) v_\lambda, \vac) = s_\lambda(1,e^{iu},e^{2iu},\dots)\,,
\end{equation}
where $\Gamma_+(u)$ is the following operator
\begin{equation}
\Gamma_+(u) = \exp \left( \sum_{n>0} \frac{1}{n} \, \frac1{1-e^{iun}} \,
\alpha_n \right) \,.
\label{G+}
\end{equation}
The coefficients inside the exponential in \eqref{G+} 
are the power-sum symmetric functions in the 
variables $1, e^{iu}, e^{2iu}, \dots$, that is, 
$$
\frac1{1-e^{iun}} = \sum_{k=0}^\infty e^{iukn} \,.
$$

Using the $u\mapsto -u$ symmetry of the right side of the GMV
formula 
\eqref{GV}, we obtain: 
  \begin{multline}\label{GVo}
\left(\prod 
\binom{(a+1)\mu_i}{\mu_i} \right)
\, \Hb\left(\mu;-1,-\frac1a,\frac1{a+1};i \sqrt{a(a+1)}\, u
\right) = \\
(-iau)^{-\ell(\mu)} 
\llang
\Gamma_+(u) \, e^{-iau \F_2 + 
\frac{iu}{2} H} \, \kket{\mu}\, , 
\end{multline}
where $H$ is the energy operator.

We perform the summation \eqref{Zloc} defining
$Z^\bullet_d(\nu,u)$
using the operator formulas \eqref{oper} and \eqref{GVo} 
for the occurring Hodge integrals and 
the formula 
$$
P_d = \sum_{|\mu|=d} \,
\frac{1} {\zz(\mu)} \,
\ket{\mu} \bra{\mu}\,,
$$
for the  orthogonal projection $P_d$ onto the space of 
vectors of energy $d$. 
The steps here exactly follow Section 3.1 of \cite{P1}.

%Here $\Pv$ is the orthogonal 
%projection onto the vacuum vector. 

The operator $P_d$ commutes with 
the operator $\F_2$. It follows that the terms
$\exp(\pm aiu\F_2)$  of \eqref{oper} and 
\eqref{GVo} cancel.  We
obtain the following all degree generating 
function in the auxiliary degree variable $Q$,  
\begin{multline}
  \label{Zoper}
\sum_d Q^d  Z^\bullet_d(\nu;u) = 
(iau)^{-n} \left[\prod z_i^{\nu_i+1}\right] \\
\lang \prod \A(z_i,iauz_i) \, e^{\alpha_1} 
\left(-\frac{Q e^{\frac{iu}2}}{iau}\right)^H \Gamma_-(u)
\rang\,. 
\end{multline}
Here, 
\begin{align*}
\Gamma_-(u) &= \Gamma_+(u)^* 
=\exp \left( \sum_{n>0} \frac{1}{n} \, \frac1{1-e^{iun}} \,
\alpha_{-n} \right) \,.
\label{G-}
\end{align*}
is the transpose of $\Gamma_+(u)$.

\subsection{Extracting the connected part}

\ssubsubsection
We must extract the degree $d$ connected part from 
the matrix element \eqref{Zoper}. 
By definition, 
$\Gamma_-(u)\, \vac$
is a linear combination of vectors of the 
form \eqref{mu}. We can also expand the 
vector 
\begin{equation}
e^{\alpha_{-1}}\, \prod \A(z_i,iauz_i)^*\, \vac
\label{prAvac}
\end{equation}
in the basis \eqref{mu}. The matrix element  \eqref{Zoper}
is then obtained via
the canonical pairing
\begin{equation}
  \label{canpar}
  \bbra{\mu} \lambda \rrang = \zz(\mu) \, \delta_{\lambda,\mu} \,. 
\end{equation}
The pairing
\eqref{canpar} may be interpreted as  an enumeration of branched
coverings of $\proj^1$ ramified over two points. Of course, a
cover of $\proj^1$ is connected only if the partition 
$$\mu=\lambda$$ has exactly one 
part. 
We will see the connected part of \eqref{Zoper} corresponds
to the contributions of the connected  covers of $\proj^1$.

\ssubsubsection
We introduce 
a weight filtration on operators on $\LV$, or,
more precisely, on the universal enveloping algebra of the 
Lie algebra $\gli$ which acts on $\LV$. 
Define the weights of the spanning elements $H^k \alpha_m$
by:
\begin{equation}
\wt H^k \alpha_m  = m + k -1  \,, \quad k=0,1,\dots, \quad m\in \Z\,,
\label{filtr}
\end{equation}
where the multiplication is in the associative algebra $\End(\infty)$.
Since
\begin{equation}
  \label{commf}
  \left[ H^a \alpha_b, H^c \alpha_d\right] = 
(bc-ad) \, H^{a+c-1} \alpha_{b+d} \,,
\end{equation}
we obtain  a well-defined filtration on the Lie algebra $\gli$ and
the associated universal enveloping algebra. 

To simplify notation, let the coefficients of $\A(z,iauz)$ be denoted by $\A_k$,
\begin{equation}
  \label{defAk}
  \A_k = [z^k]\, \A(z,iauz) \,.
\end{equation}
The weight, 
\begin{equation}
  \label{degcA}
  \wt\, \A_k  = k-1 \,,
\end{equation}
is obtained directly from the definition of $\A(z,iauz)$. 
Similarly,
\begin{equation}
  \label{smallwt}
   \A_k =  [z^{k}] \,
\A(0,iauz) +
\dots\,,
\end{equation}
where dots stand for terms of weight smaller than $k-1$. 

\ssubsubsection
The Lie algebra $\gli$ acts on $\LV$ and the filtration \eqref{filtr}
is compatible with the following filtration on $\LV$. 
We set, by definition, 
\begin{equation}
  \label{filtr2}
  \wt \ket{\mu} = |\mu| - \ell(\mu) \,,
\end{equation}
Equation \eqref{commf} implies that indeed 
\begin{equation}
  \label{degXmu}
   \wt \, X^*\!\ket{\mu} \le \wt X + \wt \ket{\mu} \,. 
\end{equation}
{}From the definitions, we have
\begin{equation}
  \label{wtvac}
   \A_k^* \, \vac = \frac{(iau)^{k}}{k!} \, 
\ket{k} + \dots\,,
\end{equation}
where dots stand for term of energy at most $k-1$ and, hence, of
weight at most $k-2$. More generally, from the commutation 
relation 
\begin{equation}
  \label{commr}
  \left[\alpha_n,\cE_{m}(s)\right] = \cs(ns) \, \cE_{m+n}(s)\,,
\end{equation}
where
\begin{equation}
  \label{defcs}
  \cs(x) = e^{x/2} - e^{-x/2} \,,
\end{equation}
we obtain 
\begin{multline}
  \label{wtcommvac}
  \left[\alpha_{\mu_1},\left[\alpha_{\mu_2},\dots \left[\alpha_{\mu_r},
\A_k\right]\right]\right]^* \, \vac =\\
 \frac{(iau)^{k}}{(k-r)!} \prod_{j=1}^r \mu_j
\ket{k-r+\sum {\mu_i}}+\dots \,,
\end{multline}
where the dots stand for terms of lower energy and weight. In particular,
the leading term vanishes if $r > k$. 

\ssubsubsection
We now apply the operators $\A(z_i,iauz_i)^*$ in 
\eqref{prAvac} in order. After each application, we expand the 
result in the basis $\ket{\mu}$. The action of the 
operator $\A_k^*$ on the vector $\ket{\mu}$ is obtained
by summing over all subsets of 
the $\mu_i$'s with which $\A_k^*$ { interacts}. 
Here, {\em interaction} denotes commutation before application
 to the vacuum. Such 
an interaction is described by equation \eqref{wtcommvac}. 
Interaction histories
can be recorded as diagrams.
As usual, extracting the connected part means
extracting the contribution of the connected diagrams. 
Since the operator $e^{\alpha_{-1}}$ 
does not interact with the operators $\A_k^*$, the operator 
$e^{\alpha_{-1}}$
 is a spectator in the extraction of the connected part. 

 From equation \eqref{wtcommvac}, we find the 
leading contribution of a connected diagram 
to the expansion of $\prod \A_{\nu_i+1}^* \, \vac$
is a constant multiple of the vector $\ket{1+|\nu|}$. 
In particular, to obtain a connected  
degree $d$ contribution, the size of the 
partition $\nu$ must be at least $d-1$. In other words,
the degree $d$ connected part of \eqref{Zoper}
vanishes if $|\nu| < d-1$.

\ssubsubsection
For the computation in the first nonvanishing case, 
\begin{equation}
  \label{snu}
  |\nu|=d-1 \,,
\end{equation}
we use formula \eqref{smallwt}.
By equation (3.12) in \cite{P1}, 
\begin{equation}
  \label{Aconj}
   \A(0,iauz_i) = e^{\alpha_1} \, \cE(iauz) \, e^{-\alpha_1} \,.
\end{equation}
Clearly,
\begin{equation}
  \label{coeffA}
  \left[z^{\nu_i+1}\right] \, \cE(iauz)
= (iau)^{\nu_i+1} \, \frac{\cP_{\nu_i+1}}
{(\nu_i+1)!}
\end{equation}
After substituting \eqref{coeffA}
 into \eqref{FZ} and \eqref{Zoper} and simplifying,
we obtain
\begin{equation}
  \label{Ioper}
  I_d(\nu;u) = \frac{(-1)^{d+1}}{2du \sin \frac{du}{2}} 
\llang e^{\alpha_1} \prod \, \frac{\cP_{\nu_i+1}}
{(\nu_i+1)!} \, \kket{d}^\circ\,,
\end{equation}
where the superscript $\circ$ denotes the connected part of the matrix
element.

The following Lemma finishes the computation and completes the proof
of Theorem \ref{GV}. 
\begin{Lemma} We have 
\begin{equation}
  \label{Hur0}
\llang e^{\alpha_1} \prod \, \frac{\cP_{\nu_i+1}}
{(\nu_i+1)!} \, \kket{d}^\circ  = 
\frac{d^{n-1}}{\prod \nu_i !} \,.
\end{equation}
\end{Lemma}

\begin{proof}
This matrix element is a relative degree $d$ Gromov--Witten 
invariant of $\proj^1$, or equivalently, a Hurwitz number with 
completed cycles insertions \cite{CC}. By condition \eqref{snu},
the invariant has genus zero and, therefore, also equals, up to a factor,
the corresponding ordinary Hurwitz number. The value of the
the genus zero Hurwitz number is well-known, see for example 
\cite{B}. Alternatively, the matrix element can be easily computed by using the commutation  
relations among the operators involved.
\end{proof}

\section{Operator formula for Hodge integrals}\label{s2}

\subsection{Operators $\bA(z,a)$}

\ssubsubsection
We begin with an operator form of the GMV formula
equivalent to equation \eqref{GVo}: 
  \begin{multline}\label{GVo2}
\left(\prod 
\binom{(a+1)\mu_i}{\mu_i} \right)
\, \Hb\left(\mu;-1,-\frac1a,\frac1{a+1}; \sqrt{a(a+1)}\, u
\right) = \\
e^{-u|\mu|/2}\,
(au)^{-\ell(\mu)} 
\llang
\Gamma_+(iu) \, e^{au \F_2} \, \kket{\mu} \,. 
\end{multline}
We will transform the
above formula by 
commuting the operators $\Gamma_+(iu)$ and  $e^{au \F_2}$, which 
fix the vacuum vector,  through the 
operators $\alpha_{-\mu_i}$.
Our strategy here follows
Section 2.2 of \cite{P1}.

\ssubsubsection 
The first conjugation 
\begin{equation}
  \label{conj1}
  e^{au \F_2} \, \alpha_{-m} \, e^{-au \F_2} = \cE_{-m}(aum)
\end{equation}
follows easily from definitions, see
Section 2.2.2 of \cite{P1}. 

The computation
of the operator
\begin{equation}
   \Gamma_+(iu) \,  \cE_{-m}(aum) \,  \Gamma_+(iu)^{-1}
\end{equation}
requires more work. We have
\begin{equation}
\Gamma_+(iu) = \prod_{n>0} \exp \left(\frac{1}{n} \, \frac1{1-e^{-un}} \,
\alpha_n \right) \,,
\label{G+2}
\end{equation}
where the factors commute. After exponentiating the relation \eqref{commr},
we obtain 
\begin{multline}
  \label{adaln}
   \exp \left(\frac{1}{n} \, \frac1{1-e^{-un}} \,
\alpha_n \right) \, \cE_{-m}(aum) \, \exp \left(-\frac{1}{n} \, \frac1{1-e^{-un}} \,
\alpha_n \right) = \\
\sum_{k\ge 0} 
\frac1{k! \,n^k} \, \left(\frac{\cs(aunm)}{1-e^{-un}}\right)^k  \cE_{-m+kn}(aum) \,.
\end{multline}
We find, 
\begin{multline}\label{bigconj}
  \Gamma_+(iu) \,  \cE_{-m}(aum) \,  \Gamma_+(iu)^{-1} = \\
\sum_{k\ge 0} \cE_{-m+k}(aum) \sum_{k_1+2k_2+3k_3+\dots=k} \frac1{\prod k_n! \, n^{k_n}} \,
\prod \left(
\frac{\cs(aumn)}{1-e^{-un}}\right)^{k_n} \,.
\end{multline}

\ssubsubsection
Let $p_k$ and $h_k$ denote the power sum and complete homogeneous
symmetric functions.
Let $\Phi$ be the specialization of the algebra of the symmetric 
functions defined by:
$$
\Phi(p_n)  = \frac{e^{aumn/2}-e^{-aumn/2}}{1-e^{-un}}\,, \quad n=1,2,\dots\, .
$$
The inner sum in the right-hand side of \eqref{bigconj} equals $\Phi(h_k)$
by standard results in the theory of symmetric functions, see \cite{Ma}, formula ($2.14'$).
Moreover, the 
equation,
$$
\Phi(h_k) = \prod_{j=1}^k \frac{e^{aum/2}-e^{-aum/2-(j-1)u}}{1-e^{-uj}} \,,
$$
is a restatement of
the $q$--binomial theorem, see \cite{Ma}, Example I.2.5.
We conclude
\begin{multline}\label{bigconj2}
  \Gamma_+(iu) \,  \cE_{-m}(aum) \,  \Gamma_+(iu)^{-1} = \\
\sum_{k\ge 0}
\left(\prod_{j=1}^k \frac{e^{aum/2}-e^{-aum/2-(j-1)u}}{1-e^{-uj}}\right)
 \cE_{-m+k}(aum)\,.
\end{multline}

\ssubsubsection
For a natural number $m$, introduce the function 
\begin{equation}
  \label{defcR}
  \cR(m,a,u) = \prod_{j=1}^m \frac{\cS((am+j-1)u)}{\cS(ju)} \,,
\end{equation}
where
$$\cS(z)=
\frac{\sinh{z/2}}{z/2} \,.$$
The definition can be extended to nonintegral values of $m$ by the 
following absolutely converging infinite product 
\begin{equation}
  \label{defcR2}
  \cR(z,a,u) = \prod_{j=1}^\infty \frac{\cS((az+j-1)u)\, \cS((j+z)u)}
{\cS((az+z+j-1)u) \,\cS(ju)\, } \,.
\end{equation}
A series expansion of this function will be discussed below in Section 
\ref{secsR}.

\ssubsubsection
Introduce the operator
\begin{multline}
  \label{defA}
  \bA(z;a) = \frac{1}{(a+1)\, u} \, \cR(z,a,u) \times \\
 \sum_{l\in\Z} 
\left(\prod_{j=1}^l \frac{e^{auz/2}-e^{-auz/2-(z+j-1)u}}{1-e^{-u(z+j)}}\right)
 \cE_{l}(auz)\,.
\end{multline}
Taking into account all prefactors, equation \eqref{bigconj2}
can be recast in the following form.

\begin{Theorem}\label{thm2}
  For positive integral values of the variables $z_i$, we have
  \begin{equation}
    \label{thrm2}
    \Hb\left(z;-1,-\frac1a,\frac1{a+1}; \sqrt{a(a+1)}\, u
\right) = 
\lang
\prod \bA(z_i;a) \rang 
  \end{equation}
\end{Theorem}

%  Since both sides of the equation make perfect sense as
%formal power series in $u$,  
%equality \eqref{thrm2} is expected to hold in 
%full generality. 
After suitable interpretation,
we expect equality \eqref{thrm2} 
to hold for all values of the variables $z_i$.
An approach along the lines of \cite{P1}
would involve establishing the commutation relations of 
the operators $\bA(z;a)$. We plan to address the topic
 in the future.

The right side of \eqref{thrm2} is less
symmetric than the left side. For example, the symmetry
with respect to 
$$
a\mapsto -a -1
$$
is not obvious from the operator formula.

\subsection{Series expansion of the function $\cR$} \label{secsR}

\ssubsubsection
Recall,
\begin{equation}
  \label{lnS}
  \ln \cS(x) = \sum_{k>0} \frac{B_{2k}}{2k\, (2k)!} \, x^{2k}\,,
\end{equation}
where $B_{m}$ are the Bernoulli numbers defined by
$$
\frac{x}{e^x-1} = \sum_{m\geq 0} \frac{B_m}{m!} \, x^m \,.
$$
{}From \eqref{lnS}, we have
\begin{equation}
  \label{lnR}
  \ln \cR(m,a,u) = \sum_{k>0} \frac{B_{2k}\, u^{2k}}{2k\, (2k)!} \sum_{j=1}^m 
\left[(am+j-1)^{2k} - j^{2k}\right] \,.
\end{equation}
The inner sum in \eqref{lnR} can be in turn computed in terms of the 
Bernoulli numbers. We obtain the following result.

\begin{Proposition}
  We have
  \begin{multline}
    \label{lnRz}
    \ln \cR(z,a,u) =\\ \sum_{0 \le l \le 2k+1} 
\frac{B_{2k} \, B_{2k-l+1} \, u^{2k}}{(2k) \, l! \, (2k-l+1)!} 
\left[(az+z)^l-(az)^l-(z+1)^l+1\right] \,.
  \end{multline}
\end{Proposition}

\ssubsubsection
Theorem \ref{thm2} implies a formula for
the 1--point function, 
\begin{equation}
  \label{1ptf}
  \Hb\left(z_1;-1,-\frac1a,\frac1{a+1}; \sqrt{a(a+1)}\, u
\right) = \frac1{(a+1)u} \, \frac{\cR(z_1,a,u)}{\cs(az_1u)} \, ,
\end{equation}
since only the constant term of the operator $\bA(z_1,u)$ contributes
to the vacuum expectation $\lang \bA(z_1,u) \rang$. 

%The series \eqref{lnRz} provides an efficient method for  computing the
%1--point function.

\end{document}